\documentclass[mathpazo]{csam}

\usepackage[svgnames, table]{xcolor}
\usepackage{booktabs}
\usepackage{multirow}
\usepackage{algorithmic}
\usepackage{algorithm}
\usepackage{amsmath}
\usepackage{subfigure}

\DeclareMathOperator*{\argmax}{argmax}

\begin{document}

\title{Compressing MIMO Channel Submatrices with Tucker Decomposition: Enabling Efficient Storage and Reducing SINR Computation Overhead}

\author[Y. Zhang et~al.]{Yuanwei Zhang\affil{1},
      Ya-Nan Zhu\affil{1}~and Xiaoqun Zhang\affil{2}\comma\corrauth}
\address{\affilnum{1}\ School of Mathematical Sciences, Shanghai Jiao Tong University, Shanghai 200240, P.R. China. \\
          \affilnum{2}\ School of Mathematical Sciences, MOE-LSC, Institute of Natural Sciences, Shanghai Jiao Tong University, Shanghai 200240, P.R. China.}
\email{{\tt xqzhang@sjtu.edu.cn} (Xiaoqun Zhang)}

\begin{abstract}
Massive multiple-input multiple-output (MIMO) systems employ a large number of antennas to achieve gains in capacity, spectral efficiency, and energy efficiency. However, the large antenna array also incurs substantial storage and computational costs. This paper proposes a novel data compression framework for massive MIMO channel matrices based on tensor Tucker decomposition. To address the substantial storage and computational burdens of massive MIMO systems, we formulate the high-dimensional channel matrices as tensors and propose a novel groupwise Tucker decomposition model. This model efficiently compresses the tensorial channel representations while reducing SINR estimation overhead. We develop an alternating update algorithm and HOSVD-based initialization to compute the core tensors and factor matrices. Extensive simulations demonstrate significant channel storage savings with minimal SINR approximation errors. By exploiting tensor techniques, our approach balances channel compression against SINR computation complexity, providing an efficient means to simultaneously address the storage and computational challenges of massive MIMO.
\end{abstract}

\maketitle

\section{Introduction}
Multiple-input and multiple-output (MIMO) technology utilizes multiple transmission and receiving antennas to exploit multipath propagation \cite{bolcskei2006mimo,jensen2004review}. It has served as the foundation for wireless and mobile networks such as the fourth (4G) and fifth (5G) generations \cite{clerckx2013mimo,kansal2019multiuser,prasad2017energy}. Compared with MIMO system, massive MIMO significantly improve the spectral and transmit power efficiency by equipping hundreds or even thousands of antennas for base station (BS). With a large number of antennas, BSs in massive MIMO system can serve multiples users simultaneously at very low signal-to-interference noise level \cite{marzetta2015massive,larsson2014massive,chataut2020massive}. For instance, 
with a sufficiently large number of antennas, linear precoding methods can achieve performance comparable to that of optimal nonlinear schemes \cite{marzetta2010noncooperative}. Moreover, if the number of BS antennas tends to be infinity, the impact of noise and intra-interference will vanish \cite{marzetta2010noncooperative}. However, the integration of a large number of antennas in a massive MIMO system also presents unprecedented challenges, particularly in terms of computation and storage.

High computational cost is incurred by signal transmission operations involving large channel matrices. Particularly, numerous operations require computational complexity with a cubic dependence on the channel size \cite{tse2005fundamentals, zeng2012linear, anand2015mode}. For example, in this paper we consider the linear minimum mean squared error (MMSE) equalization and singular value decomposition (SVD) precoding method, for achieving an optimal Signal-to-Interference Noise Ratio (SINR) among all linear schemes \cite{xie1990family}. This equalization method involves matrix-matrix multiplication and matrix inversion while the precoding method relies on SVD of channel matrices. These operations impose a substantial computational cost on massive MIMO implementations due to the large number of antennas and served users. On the other hand, channel state information (CSI) matrices are necessary for both precoding and equalization. As the number of antennas is significantly larger in massive MIMO than that in traditional MIMO systems, the CSI matrices are considerably larger. As a results, the massive MIMO system requires much higher memory capacity than conventional MIMO systems, often exceeding the size by over 100 times in practical scenarios \cite{liu2017reducing}.

Current research usually treats data compression and computational complexity reduction as two distinct objectives. Numerous studies focused on developing efficient and low-complexity algorithms for precoding/equalization. In \cite{zarei2013low,benzin2019low,thurpati2022performance}, truncation techniques were employed to approximate matrix inversion using a Taylor series expansion for achieving complexity reduction. \cite{honig2001performance} relies on the low-rankness or sparse properties inherent from the system to effectively reduce complexity. Regarding data compression methods, approaches in \cite{han2014projection, sim2016compressed} focused on reducing the size of channel data by converting it into sparse matrices. \cite{kuo2012compressive} employs dimensional reduction or compressive sensing techniques for further compression. \cite{lee2015antenna} introduces a technique that entails grouping analogous channel data from antennas and substituting them with average values and representative group patterns.

In this paper, we consider the orthogonal frequency-division multiplexing (OFDM) based massive MIMO systems. Within each channel, the presence of both line-of-sight (LOS) and none-line-of-sight (NLOS) radio waves can be represented as channel submatrices that contribute to channel matrix through linear combinations \cite{tse2005fundamentals, 3gpp.38.901}. Thus the submatrices within each channel naturally form a $3$-order tensor. Research on different tensor decomposition forms has a long history and is continuously evolving. Two well-known tensor decompositions are CANDECOMP/PARAFAC (CP) decomposition \cite{carroll1970analysis,harshman1970foundations} and Tucker decomposition \cite{tucker1966some}. CP decomposition represents a high-order tensor as a sum of rank-one tensors, while Tucker decomposition can be seen as a higher-order extension of principal component analysis (PCA) \cite{jolliffe2002principal, pearson1901liii}. Another notable decomposition is the Tensor-Train decomposition, which has gained popularity in machine learning due to its efficient implementation of basic operations \cite{oseledets2011tensor, yang2017tensor}. In the context of our problem, since signal processing and SINR calculation predominantly entail matrix product operations, we consider Tucker decomposition to be the most suitable approach for the following signal processing task, considering the benefits of reduced computational complexity when using the compressed data.

We design a Tucker decomposition based model for compressing the channel submatrices in the MIMO system. The outcome of the factorized form can not only reduce the storage cost but also reduce the computational cost of precoding, equalizer and SINR. Through numerical experiments, we demonstrate that our proposed model can achieve significant speedup of factor $6.19$ and compression ratio of $6.16$, respectively, while maintaining an acceptable level of compressed error around $10\%$. These results show the potential benefits and feasibility of our approach in practical MIMO system design.

To the best of our knowledge, this paper is the first work that tackle the compression problem of the channel submatrices to reduce both computation and memory cost. 
The advantages of our compression model can be summarized as follows:   The compressed data structure enables direct calculation of the precoder, equalizer, and SINR without reconstructing the channel matrix. This provides significant computational savings, making real-time performance feasible on devices with limited capabilities. Concurrently, compression results demonstrate substantially reduced storage requirements for channel data, enabling feasible storage on devices with constrained capacity.

The remaining sections of this paper are organized as follows. In Section {\ref{section: MCMU-MIMO System}}, we provide an introduction to the multi-cell multi-user massive MIMO system, highlighting the two main challenges addressed in this paper. We also introduce the tensor notations and explain how the Tucker decomposition method is utilized for compression. The main results of this study are presented in Section {\ref{section: Proposed Methods}}. Firstly, we introduce a simple compression model based on Tucker decomposition and discuss the key ideas for achieving faster computations. Then, we propose the Groupwise Tucker  compression model and present the corresponding algorithms for compressing channel data in Section {\ref{section: MT model}}. Additionally, we provide a detailed analysis of the reduction in computation complexity of SINR in Section {\ref{section: SINR Complexity Analysis}}. In Section {\ref{section: Numerical Results}}, we present the numerical results obtained from various experiments conducted to evaluate the performance of the proposed models. Finally, we conclude the paper in Section {\ref{section: Conclusion and Discussion}}, summarizing the key findings and discussing the implications of our work.

\section{Background}{\label{section: MCMU-MIMO System}}
In this section, we first present some related notations and introduce the background for MIMO system.
\subsection{Notations}
\begin{itemize}
    \item Matrices are denoted by bold capital letters (e.g., $\boldsymbol{H}$), vectors are denoted by bold lowercase letters (e.g., $\boldsymbol{x}$) and tensors are denoted by bold calligraphic letters (e.g., $\boldsymbol{\mathcal{X}}$). The index set $[N]$ denotes the set $\{1,2,\dots,N\}$, where $N$ is a positive integer. $\sqrt{-1}$ represents the imaginary unit. 
    \item $\boldsymbol{A}^*$ represents the hermitian matrix of $\boldsymbol{A}$. $[\boldsymbol{A}]_r$ represents the $r$-th column of matrix $\boldsymbol{A}$ and $[\boldsymbol{A}]_{1:r}$ represents the first $r$ columns of $\boldsymbol{A}$. $\boldsymbol{I}_N$ denotes the identity matrix with size $N\times N$.
    \item The bold symbols $\boldsymbol{H}^{[i, j]}_k\in \mathbb{C}^{M\times N}$ and $\boldsymbol{\boldsymbol{\mathcal{X}}}_k^{[i,j]}\in \mathbb{C}^{M\times N\times P}$ represent the channel matrix and tensor composed of channel submatrices, respectively. In this context, $M$ and $N$ refer to the number of antennas for the receiver and transmitter, respectively, while $P$ represents the number of submatrices. The indices $[i,j]$ correspond to specific users and $k$ serves as the index for the base station. Similarly, we have $\boldsymbol{\widetilde{H}}^{[i,j]}_k\in \mathbb{C}^{m\times n}$ and  $\boldsymbol{\boldsymbol{\mathcal{G}}}^{[i,j]}_k\in \mathbb{C}^{m\times n\times p}$, which represent the compressed channel matrix and channel tensor ($m\leq M, n\leq N, p \leq P$). 
\end{itemize}

\subsection{Signal Transmission and SINR}
\begin{figure}[h!]
    \centering
    \includegraphics[width = .7\textwidth]{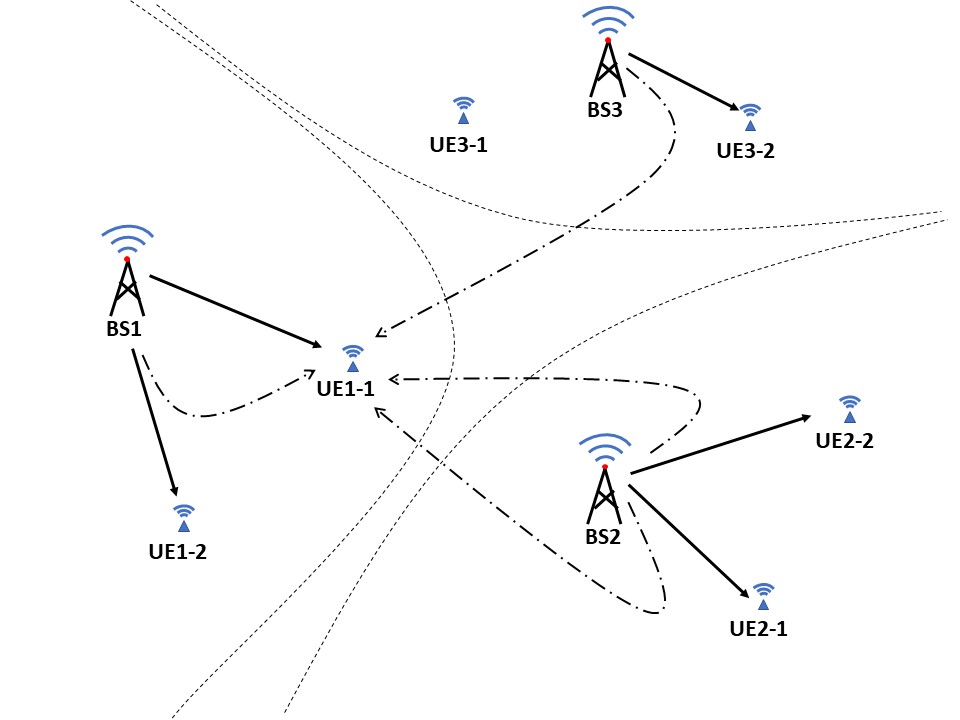}
    \caption{Illustration of a MIMO system. In this MIMO system, each cellular network cell is equipped with a single base station serving two users. The solid lines denote the desired signal paths and the dashed lines signify interference. Specifically, UE1-1 receives the desired signal from BS1, experiencing both intra-cell interference from BS1 and co-channel interference from BS2 and BS3.}
    \label{Illustration: MIMO system}
\end{figure}

{\label{section: Signal transmit and SINR calculation}}
As depicted in Figure {\ref{Illustration: MIMO system}}, we consider a  general multi-cell multi-user MIMO system, comprising $J$ Base Stations (BSs), with one BS allocated per cell. Let $K_i$ be the number of scheduled users for $i$th-BS. Each link consists of a BS equipped with $N$ transmit antennas and a user equipment (UE) equipped with $M$ receive antennas. Here, we denote the $M\times N$ channel matrix from $k$-th BS to $j$-th UE in $i$-th cell as $\boldsymbol{H}^{[i,j]}_k,\ i,k \in [J]$ and $j \in [K_i]$.

For universality, we consider both intra-cell and inter-cell interference. Intra-cell interference occurs when multiple UEs are concurrently scheduled using spatial division technique within the same cell, while the inter-cell interference arises from BSs in neighboring cells when they serve their UEs.
Let $L$ denote the spatial streams number, then the received signal $\boldsymbol{y}^{[i,j]}$ for $j$-th UE in $i$-th cell can be expressed as the sum of four distinct components:

\begin{equation}
    \label{equ: signal receiving}
    \boldsymbol{y}^{[i,j]} = \underbrace{\boldsymbol{H}_{i}^{[i,j]}\boldsymbol{P}^{[i,j]}\boldsymbol{x}^{[i,j]}}_{\text{desired signal}} + \underbrace{\sum_{l = 1, l\neq i}^{K_i}\boldsymbol{H}^{[i,j]}_i \boldsymbol{P}^{[i,l]}\boldsymbol{x}^{[i,l]}}_{\text{intra-cell interference}} + \underbrace{\sum_{k = 1,k\neq i}^J \sum_{l = 1}^{K_k} \boldsymbol{H}_{k}^{[i,j]} \boldsymbol{P}^{[k,l]} \boldsymbol{x}^{[ k,l]}}_{\text{inter-cell interference}} + \underbrace{\boldsymbol{z}^{[i,j]}}_{\text{noise}}, 
\end{equation}
where $\boldsymbol{P}^{[i,j]}$ is the $N\times L$ precoding matrix and $\boldsymbol{x}^{[i,j]}$ is the $L \times 1$ transmit-signal vector for $j$-th UE in $i$-th cell. $\boldsymbol{z}^{[i,j]}$ is the additive white Gaussian noise.

For each $\boldsymbol{H}^{[i,j]}_k,\ i,k \in [J]$ and $j \in [K_i]$, the  channel matrix is generated from submatices with different time-delay.  Denote the channel matrix with time-varying as $\boldsymbol{H}(t)$, then $\boldsymbol{H}(t)$ is comprised of a linear combination of  a line-of-sight (LOS) submatrix $\boldsymbol{H}^{LOS}(t)$ and multiple non-line-of-sight (NLOS) submatrices $\boldsymbol{H}^{NLOS}_l(t), l = 1,\dots, P-1$, as follows:
\begin{equation}
    \boldsymbol{H}(t) = c_0 (t-t_0) \boldsymbol{H}^{LOS}(t_0) + \sum_{l=1}^{P-1} c_l(t - t_0) \boldsymbol{H}^{NLOS}_{l}(t_0),
    \label{equ: channel submatrices}
\end{equation}
where $P$ is the number of submatrices and $c_l(t)$ is the time-varying coefficients.

For the receiver, once it has received the signal $\boldsymbol{y}^{[i,j]}$, the UE utilizes a receive filter matrix $\boldsymbol{W}^{[i,j]}$ to reconstruct the desired transmit-signal vector $\boldsymbol{x}^{[i,j]}$, i.e., the estimated desired signal is computed by
$\hat{\boldsymbol{x}}^{[i,j]} = (\boldsymbol{W}^{[i,j]})^* \boldsymbol{y}^{[i,j]}$. 
Here we consider MMSE equalizer \cite{jiang2011performance} and $\boldsymbol{W}^{[i,j]}$ is defined as follows:

\begin{equation}
    \label{equ: MMSE receive filter}
    \boldsymbol{W}^{[i,j]} = (\boldsymbol{Q}^{[i,j]})^{-1} \boldsymbol{H}_{i}^{[i,j]}\boldsymbol{P}^{[i,j]},
\end{equation}
where $\boldsymbol{Q}^{[i,j]}$ is the total received signal covariance matrix of $j$-th UE in $i$-th cell. Assuming the standard variation of the Gaussian noise as $\sigma$,  $\boldsymbol{Q}^{[i,j]}$ can be expressed as:
\begin{equation}
    \boldsymbol{Q}^{[i,j]} = \sum_{k = 1}^J \sum_{l = 1}^{K_k} \boldsymbol{H}_{k}^{[i,j]} \boldsymbol{P}^{[k,l]}(\boldsymbol{P}^{[k,l]})^*(\boldsymbol{H}_{k}^{[i,j]})^* + \sigma^2\boldsymbol{I}_M. 
    \label{equ: covariance matrix}
\end{equation}

In accordance with the notations provided above, the SINR of the $r$-th spatial stream of $j$-th UE in $i$-th cell is calculated as follows:
\begin{equation}
    \label{equ: SINR formula}
        {SINR}^{[i,j]}_r = \frac{\boldsymbol{s}^{[i,j]}_r}{\left[\boldsymbol{W}^{[i,j]}\right]_r^* \boldsymbol{Q}^{[i,j]}\left[\boldsymbol{W}^{[i,j]}\right]_r-\boldsymbol{s}^{[i,j]}_r}
\end{equation}
where 
\begin{equation}\label{equ: s_r^ij}
    \boldsymbol{s}^{[i,j]}_r = \left|\left[\boldsymbol{W}^{[i,j]}\right]_r^* \boldsymbol{H}^{[i,j]}_i\left[\boldsymbol{P}^{[i,j]}\right]_r\right|^2
\end{equation}
is the signal power of $r$-th spatial stream. Here $r \in [L]$ and $\left[\cdot\right]_r$ denotes the $r$-th column vector of the matrix. The denominator represents the power of total interference corresponding to $r$-th spatial stream. 

The equations {(\ref{equ: covariance matrix})} and {(\ref{equ: SINR formula})} show that the value of SINR depends on the precoding method of transmit signal, receive filter method, and the number of spatial streams, etc. 
In this paper, we consider the SVD precoding method, which can achieve better channel capacity \cite{telatar1999capacity}. Denote $\boldsymbol{U}^{[i,j]}\boldsymbol{\Sigma}^{[i,j]} (\boldsymbol{V}^{[i,j]})^*$ is the truncated-$L$ SVD of $\boldsymbol{H}^{[i,j]}_i$, where $\boldsymbol{V}^{[i,j]}$ represents the first $L$ right singular vectors of $\boldsymbol{V}^{[i,j]}$. Then, the SVD precoding method is given by setting
\begin{equation}
\boldsymbol{P}^{[i,j]} = \boldsymbol{V}^{[i,j]}.
\label{equ: precoding matrix}
\end{equation}
{Then the filter matrix $ \boldsymbol{W}^{[i,j]} = (\boldsymbol{Q}^{[i,j]})^{-1} \boldsymbol{U}^{[i,j]}\boldsymbol{\Sigma}^{[i,j]}$ and the signal power, $\boldsymbol{s}^{[i,j]}_r = (\sigma_r^{[i,j]})^4\cdot \left| \left[\boldsymbol{U}^{[i,j]}\right]_r^* (\boldsymbol{Q}^{[i,j]})^{-1} \left[\boldsymbol{U}^{[i,j]}\right]_r \right|^2$, where $\sigma_r^{[i,j]}$ is the $r$-th singular value. It can be seen that both the calculation of the receive filter $\boldsymbol{W}^{[i,j]}$ and {SINR} require intricated operations, including the inversion of the covariance matrix $\boldsymbol{Q}^{[i,j]}$ and SVD of channel matrices.}

These above operations exhibit a substantial level of computational complexity, which becomes even more challenging in the context of uplink signal transmission as the size of the covariance matrix scales with the number of antennas at the base station. This computational challenge becomes particularly undesirable in modern applications like 5G technology, where an increasing number of antennas are utilized to mitigate growing path loss. This, serves as the primary motivation of this paper, aimed at reducing these computational complexities effectively. The second purpose of this paper is to reduce the storage overhead. In practices one need to store the submatrices in \eqref{equ: channel submatrices} of all channels at a fix time point for generating the channel matrix. However, those submatrices keeps the same size with channel matrices and have a large number, therefore, it imposes a substantial storage burden, making it impractical, particularly when dealing with an increasing number of users.

\subsection{Preliminaries on Tensors}{\label{section: preliminaries}}
Before proposing the data compression model, we first introduce some notations for tensor based operations and review some related algorithms, readers may refer to  \cite{kolda2009tensor} for more details.

\vskip 2mm
\begin{itemize}
\item \textit{Order of tensor, also known as way or mode,} refers to the number of dimensions. For example, a vector is a $1$-order tensor and a matrix is a $2$-order tensor.
\vskip 2mm
\item \textit{Inner product} of two complex $3$-order tensors $\boldsymbol{\mathcal{X}}$ and $\boldsymbol{\mathcal{Y}}$ with the same size is defined by
\begin{equation*}
    \left<\boldsymbol{\mathcal{X}}, \boldsymbol{\mathcal{Y}}\right> = \sum\limits_{i_1, i_2, i_3} x_{i_1, i_2, i_3} \bar{y}_{i_1, i_2, i_3}.
\end{equation*}
The induced tensor norm is $\|\boldsymbol{\mathcal{X}}\| = \sqrt{\left<\boldsymbol{\mathcal{X}}, \boldsymbol{\mathcal{X}}\right>}$.

\vskip 2mm
\item \textit{Matricization, or unfolding}, is the reordering of all elements of a tensor into a matrix. For example, the mode-$1$ matricizaion of $3$-order tensor $\boldsymbol{\mathcal{X}}\in \mathbb{C}^{I_1\times I_2\times I_3 }$ is denoted as $\boldsymbol{\mathcal{X}}_{(1)}\in \mathbb{C}^{I_1\times I_2 I_3}$. It is the rearrangement of vectors of length $I_1$ into a matrix according to the order. More precisely, $[\boldsymbol{\mathcal{X}}_{(1)}]_{i_1, i_2+n(i_3-1)} = \boldsymbol{\mathcal{X}}_{i_1, i_2, i_3}$.

\vskip 2mm
\item \textit{Tensor-matrix product.} Denote $\boldsymbol{\mathcal{X}}\times_1 \boldsymbol{U}$ as the $1$-mode product of the tensor $\boldsymbol{\mathcal{X}}\in \mathbb{C}^{I_1\times I_2\times I_3}$ with a matrix $\boldsymbol{U}\in \mathbb{C}^{J\times I_1}$, then $\boldsymbol{\mathcal{X}}\times_1 \boldsymbol{U}\in \mathbb{C}^{J\times I_2\times I_3}$. Elementwisely, we have
\begin{equation*}
     (\boldsymbol{\mathcal{X}}\times_{1}\boldsymbol{U})_{j,i_{2},i_{3}} = \sum\limits_{i_1 = 1}^{I_1}x_{i_1,i_2,i_3}u_{j,i_1}.
\end{equation*}
It is straightforward to verify that,
\begin{equation*}
    \boldsymbol{\mathcal{Y}}=\boldsymbol{\mathcal{X}} \times_{1} \boldsymbol{U} \quad \Leftrightarrow \quad \boldsymbol{\mathcal{Y}}_{(1)}=\boldsymbol{U} \boldsymbol{\mathcal{X}}_{(1)},
\end{equation*}
\vskip 2mm
\item \textit{Tucker rank.} Denote $R_n = \operatorname{rank}(\boldsymbol{\mathcal{X}}_{(n)})$, then the vector $\operatorname{rank_T}(\boldsymbol{\mathcal{X}}) = (R_1, R_2, \dots, R_N)$ is the Tucker rank of $N$-order tensor $\boldsymbol{\mathcal{X}}$.
\vskip 2mm
\end{itemize}

With those definitions and operations, We will first introduce the Tucker decomposition method and how directly use Tucker decomposition to compress channel tensors.

\begin{itemize}
\item \textit{Tucker decomposition}. 
Tucker decomposition is a generalization of the PCA method in higher order perspective \cite{tucker1966some}. For a $3$-order tensor $\boldsymbol{\mathcal{X}}\in \mathbb{C}^{I_1\times I_2\times I_3}$, 
it decomposed $\boldsymbol{\mathcal{X}}$ into a smaller tensor $\boldsymbol{\mathcal{G}}$ multiplied three orthogonal factor matrices $\boldsymbol{A}\in \mathbb{C}^{I_1\times R_1},\boldsymbol{B}\in\mathbb{C}^{I_2\times R_2},\boldsymbol{C}\in \mathbb{C}^{I_3\times R_3} (I_i\geq R_i, i = 1,2,3)$, i.e., $\boldsymbol{\mathcal{X}} \approx \boldsymbol{\mathcal{G}} \times_1  \boldsymbol{A} \times_2  \boldsymbol{B} \times_3  \boldsymbol{C}$ through the following optimization problem

\begin{equation}
\begin{array}{c}
\min\limits_{\boldsymbol{A}, \boldsymbol{B}, \boldsymbol{C}, \boldsymbol{\mathcal{G}}}  \left\|\boldsymbol{\mathcal{X}} - \boldsymbol{\mathcal{G}}\times_1 \boldsymbol{A}\times_2 \boldsymbol{B}\times_3 \boldsymbol{C}\right\|_{F}^2 \\
\text {s.t. }  \boldsymbol{A}^{*} \boldsymbol{A}=\boldsymbol{I}_{R_1},\boldsymbol{B}^{*} \boldsymbol{B}=\boldsymbol{I}_{R_2},\boldsymbol{C}^{*} \boldsymbol{C}=\boldsymbol{I}_{R_3}
\end{array}
\label{model: general_Tucker decomposition}
\end{equation}
There are two popular algorithms available for solving {(\ref{model: general_Tucker decomposition})}: the High Order Singular Value Decomposition (HOSVD) algorithm \cite{de2000multilinear} and the Higher Order Orthogonal Iteration (HOOI) algorithm \cite{de2000best}. The HOSVD algorithm is an extension of matrix SVD, it leverages the principal left singular vectors from each mode matricization of the target tensor as the factor matrix. 

\begin{minipage}{\linewidth}
\begin{algorithm}[H]
	\renewcommand{\algorithmicrequire}{{Input:}}
	\renewcommand{\algorithmicensure}{{Output:}}
	\caption{HOSVD Algorithm.}
	\label{alg: HOSVD Algorithm}
	\begin{algorithmic}[1]
        \REQUIRE Target $3$-order tensor $\boldsymbol{\mathcal{X}}\in \mathbb{C}^{I_1\times I_2\times I_3}$. The parameters $R_1, R_2, R_3$. 
        \STATE $\boldsymbol{A} \longleftarrow R_1 \text{ leading left singular vectors of } \boldsymbol{\mathcal{X}}_{(1)}$.
        \STATE $\boldsymbol{B} \longleftarrow R_2 \text{ leading left singular vectors of }\boldsymbol{\mathcal{X}}_{(2)}$.
        \STATE $\boldsymbol{C} \longleftarrow R_3 \text{ leading left singular vectors of }\boldsymbol{\mathcal{X}}_{(3)}$.
        \STATE $\boldsymbol{\mathcal{G}} = \boldsymbol{\mathcal{X}} \times_1 \boldsymbol{A}^*\times_2 \boldsymbol{B}^*\times_3 \boldsymbol{C}^*$
    \ENSURE Core tensor $\boldsymbol{\mathcal{G}}$, factor matrices $\boldsymbol{A}$, $\boldsymbol{B}$, $\boldsymbol{C}$.
	\end{algorithmic}
\end{algorithm}
\end{minipage}

The solutions obtained from the HOSVD algorithm are usually not optimal solution of {(\ref{model: general_Tucker decomposition})}. The HOOI algorithm incorporates the HOSVD algorithm as an initialization method and employs an alternating iteration approach to achieve a more accurate solution \cite{de2000best}.

\end{itemize}

\section{Methods}{\label{section: Proposed Methods}}

\subsection{Channel Matrices in tensor form}{\label{section: channel submatrices}}

Stacking the submatrices $\boldsymbol{H}^{LOS}(t)$ and $\boldsymbol{H}^{NLOS}_l(t), l = 1,\dots, P-1$
 in \eqref{equ: channel submatrices} along with third-mode as a tensor. \eqref{equ: channel submatrices} can be reformulated as the following tensor-vector product
\begin{equation}
\boldsymbol{H}(t) = \boldsymbol{\mathcal{X}}(t_0)\times_3 \boldsymbol{c}^*(t - t_0),
\label{equ: channel submatrices, tensor}
\end{equation}
where $\boldsymbol{c}(t) = \left[c_1(t), c_2(t), \dots, c_P(t)\right]'\in \mathbb{C}^{P\times 1}$ is the coefficient vector. Denote the tensor for each channel as $\boldsymbol{\mathcal{X}}^{[i,j]}_k$ and we want to apply decomposition method to compress those channel tensors while reduce the SINR estimation overhead,

The key idea in our proposed compression model to reducing computational complexity lies in leveraging the orthogonality of factor matrices within the decomposition form and canceling them out during SINR calculations. As a result, the sizes of matrix involved in the computations correspond to the compressed scale, thereby decreasing computational complexity. 
\begin{itemize}
    \item \textbf{Individual Tucker compression.} 
    For each 3-order tensor $\boldsymbol{\mathcal{X}}^{[i,j]}_k\in \mathbb{C}^{M\times N\times P}$. The tucker decomposition find the factor matrices $\boldsymbol{A}^{[i,j]}_k \in\mathbb{C}^{M\times m}, \boldsymbol{B}^{[i,j]}_k \in\mathbb{C}^{N\times n}, \boldsymbol{C}^{[i,j]}_k \in\mathbb{C}^{P\times p}$ and core tensor $\boldsymbol{\mathcal{G}}^{[i,j]}_k \in\mathbb{C}^{m\times n\times p}$ by solving the following problem

\begin{equation}
\begin{array}{c}
\min\limits_{\boldsymbol{A}^{[i,j]}_k, \boldsymbol{B}^{[i,j]}_k, \boldsymbol{C}^{[i,j]}_k, \boldsymbol{\mathcal{G}}^{[i,j]}_k}  \left\|\boldsymbol{\mathcal{X}}^{[i,j]}_k - \boldsymbol{\mathcal{G}}^{[i,j]}_k\times_1 \boldsymbol{A}^{[i,j]}_k\times_2 \boldsymbol{B}^{[i,j]}_k\times_3 \boldsymbol{C}^{[i,j]}_k\right\|_{F}^2 \\
\text {s.t. }  (\boldsymbol{A}^{[i,j]}_k)^{*} \boldsymbol{A}^{[i,j]}_k=\boldsymbol{I}_m,(\boldsymbol{B}^{[i,j]}_k)^{*} \boldsymbol{B}^{[i,j]}_k=\boldsymbol{I}_n,(\boldsymbol{C}^{[i,j]}_k)^{*} \boldsymbol{C}^{[i,j]}_k=\boldsymbol{I}_p
\end{array}
\label{model: Tucker decomposition}
\end{equation}

In the individually Tucker decomposition model {(\ref{model: Tucker decomposition})}, the mode-$1, 2$ factor matrices $\boldsymbol{A}^{[i,j]}_k, \boldsymbol{B}^{[i,j]}_k$ for different channels are distinct and, therefore, cannot be eliminated during SINR calculation. To address this, an straightforward approach is to make the mode-$1, 2$ factor matrices of different channels are identical. For instance, in the case of channel matrix $\boldsymbol{H}^{[i,j]}_k$ between user-$[i,j]$ and base-$k$, the compressed channel matrix can be structured in the following manner:
\begin{equation}
\begin{array}{c}
\displaystyle
\boldsymbol{H}^{[i,j]}_k \approx \boldsymbol{A} \boldsymbol{\widetilde{H}}_{k}^{[i,j]} \boldsymbol{B}^*, 
\text{ with } \boldsymbol{A}^*\boldsymbol{A} = \boldsymbol{I}_m, \boldsymbol{B}^* \boldsymbol{B} = \boldsymbol{I}_n.\end{array}
\label{equ: MT1 matrix factorization}
\end{equation}

With the utilization of this factorization form, the received signal can be conveniently reformulated. As an example, let's consider the term $\boldsymbol{H}_{k}^{[i,j]} \boldsymbol{P}^{[k,l]} \boldsymbol{x}^{[ k,l]}$ in {(\ref{equ: signal receiving})}, with the SVD precoding method $\boldsymbol{P}^{[k,l]} = \boldsymbol{V}^{[k,l]}$.
\begin{equation}
    \boldsymbol{H}_{k}^{[i,j]} \boldsymbol{P}^{[k,l]} \boldsymbol{x}^{[ k,l]}\approx\boldsymbol{A} \boldsymbol{\widetilde{H}}_{k}^{[i,j]} \boldsymbol{B}^* \boldsymbol{V}^{[k,l]} \boldsymbol{x}^{[k,l]} = \boldsymbol{A} \boldsymbol{\widetilde{H}}_{k}^{[i,j]} \boldsymbol{\widetilde{V}}^{[k,l]}\boldsymbol{x}^{[k,l]}
\end{equation}
In this context, the precoding matrix is expressed as $\boldsymbol{V}^{[k,l]} = \boldsymbol{B} \boldsymbol{\widetilde{V}}^{[k,l]}$, where $\boldsymbol{\widetilde{V}}^{[k,l]}$ represents the principal left singular vector matrix of $\boldsymbol{\widetilde{H}}_{k}^{[i,j]}$. Therefore, the factor matrix $\boldsymbol{B}$ can be eliminated during the precoding process. Similarly, the factor matrix $\boldsymbol{A}$ can be eliminated during equalization. 
As a result, SINR computation exclusively involves the compressed and reduced-scale matrix $\boldsymbol{\widetilde{H}}_{k}^{[i,j]}$. This factorization form significantly accelerates the SINR calculation process.

\item \textbf{Shared Tucker compression.}
Similar to {(\ref{model: Tucker decomposition})}, one can solve the following optimization problem to obtain the compressed channel tensors

\begin{equation}
\label{model: MT1}
\begin{array}{c}
    \min\limits_{\boldsymbol{A},\boldsymbol{B},\{\boldsymbol{C}^{[i,j]}_k\},\{\boldsymbol{\mathcal{G}}^{[i,j]}_k\}}  \sum\limits_{[i,j]}\sum\limits_{k}\left\|\boldsymbol{\mathcal{X}}^{[i,j]}_k -\boldsymbol{\mathcal{G}}^{[i,j]}_k\times_1 \boldsymbol{A} \times_2 \boldsymbol{B} \times_3 \boldsymbol{C}^{[i,j]}_k\right\|_F^2\\
    \text{s.t. }{\boldsymbol{A}}^*\boldsymbol{A} = \boldsymbol{I}_m, \boldsymbol{B}^* \boldsymbol{B} = \boldsymbol{I}_n, {\boldsymbol{C}^{[i,j]}_k}^*\boldsymbol{C}^{[i,j]}_k = \boldsymbol{I}_p, \text{ for all } i, j, k.
\end{array},
\end{equation}

Nonetheless, as the number of BSs and UEs increases, achieving a suitably minimized approximation error in {(\ref{model: MT1})} becomes challenging due to the shared mode-$1$ and mode-$2$ factor matrices $\boldsymbol{A}$ and $\boldsymbol{B}$. Consequently, attaining a sufficiently accurate compressed solution becomes difficult.

\end{itemize}

\subsection{Groupwise-Tucker Compression Model}{\label{section: MT model}}

For the above two models can not achieve a good data compression and computation acceleration within an acceptable error rate, we consider a groupwise Tucker compression method. For a given channel matrix $\boldsymbol{H}^{[i,j]}_k$ between user $[i,j]$ and base station $k$, the matrix factorization is expressed in the following form:
\begin{equation}
\boldsymbol{H}_{k}^{[i,j]}
\approx \boldsymbol{A}_{[i,j]} \boldsymbol{\widetilde{H}}_{k}^{[i,j]} {\boldsymbol{B}_k}^*, \text{ with }{\boldsymbol{A}_{[i,j]}}^* \boldsymbol{A}_{[i,j]} = \boldsymbol{I}_m,\ {\boldsymbol{B}_k}^* \boldsymbol{B}_k = \boldsymbol{I}_n,  
\label{equ: three matrix factorization}
\end{equation}
In {(\ref{equ: three matrix factorization})} and
Figure {\ref{Illustration: channel compression}}, it can be observed that the left factor matrix $\boldsymbol{A}_{[i,j]}$ represents the factor matrix for the serving and interference channels of user $[i,j]$, while the right factor matrix $\boldsymbol{B}_k$ corresponds to the factor matrix for all channels from base station $k$. Compared to the factorization form in {(\ref{equ: MT1 matrix factorization})}, this factorization form incorporates more mode-$1,2$ factor matrices, resulting in reduced approximation errors. In Section {\ref{section: SINR Complexity Analysis}}, we will demonstrate that this relaxation still preserves the property of eliminating the factor matrices ${\boldsymbol{A}_{[i,j]}}$ and ${\boldsymbol{B}_k}$ during SINR calculations.

\begin{center}
\begin{minipage}{.9\linewidth}
\begin{figure}[H]
    \centering
    \includegraphics[width = .8\textwidth]{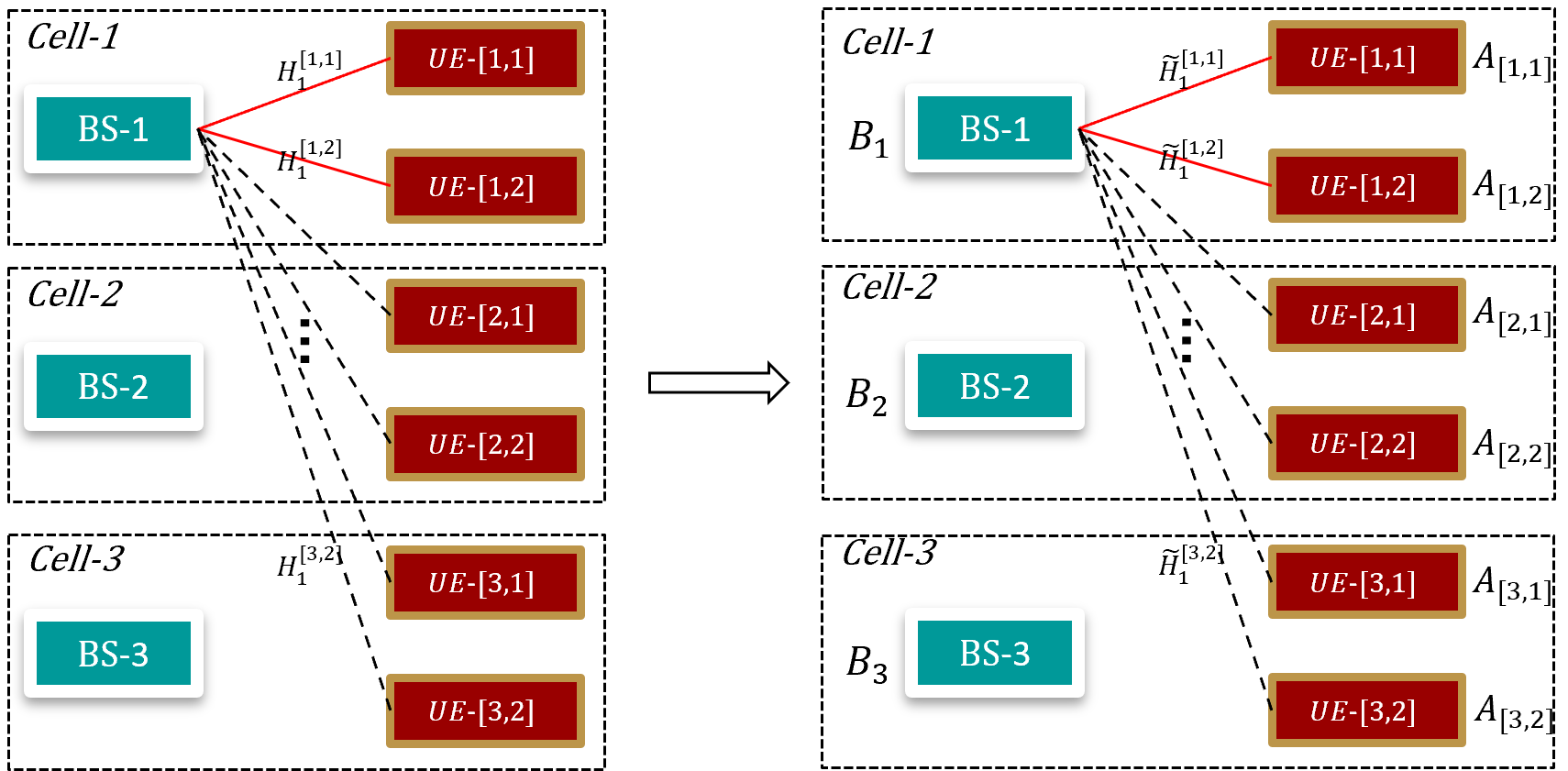}
    \caption{Channel matrices of MIMO system with 3 bases and 2 users per cell before and after compression. The serving channel matrices between BS-$1$ and UE-$[1,1]$ and UE-$[1,2]$ are represented by two red lines, while the interference channel matrices from BS-$1$ and UEs in other cells are represented by black dashed lines.}
    \label{Illustration: channel compression}
\end{figure}
\end{minipage}
\end{center}

Denote the channel tensor between user-$[i,j]$ and base-$k$ as $\boldsymbol{\mathcal{X}}^{[i,j]}_k$, and its coefficient vector in {(\ref{equ: channel submatrices})} as $\boldsymbol{c}_k^{[i,j]}$, then the Groupwise Tucker compression model factorizes $\boldsymbol{\mathcal{X}}^{[i,j]}_k$ in the following form
\begin{equation}
\begin{array}{c}
\displaystyle
\boldsymbol{\mathcal{X}}^{[i,j]}_k \approx \boldsymbol{\mathcal{G}}^{[i,j]}_k\times_1 \boldsymbol{A}_{[i,j]} \times_2 \boldsymbol{B}_k \times_3 \boldsymbol{C}^{[i,j]}_k\\
    \text{with } {\boldsymbol{A}_{[i,j]}}^*\boldsymbol{A}_{[i,j]} = \boldsymbol{I}_m, {\boldsymbol{B}_k}^* \boldsymbol{B}_k = \boldsymbol{I}_n, {\boldsymbol{C}^{[i,j]}_k}^*\boldsymbol{C}^{[i,j]}_k = \boldsymbol{I}_p.\end{array}
\label{equ: tensor approximation for SINR reduction}
\end{equation}
Plug the {(\ref{equ: tensor approximation for SINR reduction})} into   {(\ref{equ: channel submatrices, tensor})} yields
\begin{equation}
    \label{equ: compressed matrix SINR reduction}
    \boldsymbol{H}^{[i,j]}_k = \boldsymbol{\mathcal{X}}^{[i,j]}_k \times_3 {\boldsymbol{c}_k^{[i,j]}}^* \approx  \boldsymbol{A}_{[i,j]} (\boldsymbol{\mathcal{G}}^{[i,j]}_k\times_3 ({\boldsymbol{c}_k^{[i,j]}}^* \boldsymbol{C}^{[i,j]}_k)) \boldsymbol{B}_k^*.
\end{equation}
The compressed channel matrix in {(\ref{equ: three matrix factorization})} can be represented as $\boldsymbol{\widetilde{H}}^{[i,j]}_k = \boldsymbol{\mathcal{G}}^{[i,j]}_k \times_3 ({\boldsymbol{c}_k^{[i,j]}}^* \boldsymbol{C}^{[i,j]}_k)$. Comparing this with {(\ref{equ: channel submatrices, tensor})}, it appears that an additional matrix-vector multiplication is required to obtain the compressed channel matrix. However, since both $\boldsymbol{\mathcal{G}}^{[i,j]}_k$ and $\boldsymbol{C}^{[i,j]}_k$ are compressed, the overall computational cost is lower than the tensor-vector multiplication in {(\ref{equ: channel submatrices, tensor})}.

\vskip 2mm

To get the factor matrices and corresponding core tensor, analogous to Tucker decomposition, we define $f_k^{[i,j]} := \left\|\boldsymbol{\mathcal{X}}^{[i,j]}_k -\boldsymbol{\mathcal{G}}^{[i,j]}_k\times_1 \boldsymbol{A}_{[i,j]} \times_2 \boldsymbol{B}_k \times_3 \boldsymbol{C}^{[i,j]}_k\right\|_F^2$ as the approximation error of channel tensor $\boldsymbol{\mathcal{X}}^{[i,j]}_k$ and consider the following model:

\begin{equation}
\begin{array}{c}
\displaystyle
    \min\ f:= \sum\limits_{[i,j]}\sum\limits_{k}f_k^{[i,j]} = \sum\limits_{[i,j]}\sum\limits_{k}\left\|\boldsymbol{\mathcal{X}}^{[i,j]}_k -\boldsymbol{\mathcal{G}}^{[i,j]}_k\times_1 \boldsymbol{A}_{[i,j]} \times_2 \boldsymbol{B}_k \times_3 \boldsymbol{C}^{[i,j]}_k\right\|_F^2\\
    \text{s.t. }{\boldsymbol{A}_{[i,j]}}^*\boldsymbol{A}_{[i,j]} = \boldsymbol{I}_m, {\boldsymbol{B}_k}^* \boldsymbol{B}_k = \boldsymbol{I}_n, {\boldsymbol{C}^{[i,j]}_k}^*\boldsymbol{C}^{[i,j]}_k = \boldsymbol{I}_p, \text{ for all }i,j,k.
\end{array},
\label{model: MT_origin}
\end{equation}

As analyzed in \cite{kolda2009tensor}, when minimizing such an objective with orthogonal constraints on the factor matrices, it is more efficient to transform this problem into an equivalent maximization problem. It is easy to verify that if the factor matrices $\boldsymbol{A}_{[i,j]},\boldsymbol{B}_k,\boldsymbol{C}^{[i,j]}_k$ are fixed, then the optimal $\boldsymbol{\mathcal{G}}^{[i,j]}_k$ is:
\begin{equation}\label{equ: kernel}
    \boldsymbol{\mathcal{G}}^{[i,j]}_k = \boldsymbol{\mathcal{X}}^{[i,j]}_k\times_1 {\boldsymbol{A}_{[i,j]}}^* \times_2 {\boldsymbol{B}_k}^* \times_3 {\boldsymbol{C}^{[i,j]}_k}^*.
\end{equation}
By substituting  $\boldsymbol{\mathcal{G}}^{[i,j]}_k$ in   {(\ref{equ: kernel})} into $f_k^{[i,j]}$, one gets
\begin{equation} \label{equ: min-max substitite}
\begin{aligned}
    f_k^{[i,j]} =& \left\|\boldsymbol{\mathcal{X}}^{[i,j]}_k \right\|_F^2 - \left\|\boldsymbol{\mathcal{X}}^{[i,j]}_k\times_1 {\boldsymbol{A}_{[i,j]}}^* \times_2 {\boldsymbol{B}_k}^* \times_3 {\boldsymbol{C}^{[i,j]}_k}^*\right\|_F^2\\
    :=& \left\|\boldsymbol{\mathcal{X}}^{[i,j]}_k \right\|_F^2 - g_k^{[i,j]} 
\end{aligned}
\end{equation}
Thus, {(\ref{model: MT_origin})} is equivalent to the following maximization problem:
\begin{equation}
\begin{array}{c}
    \displaystyle
    \max\ g:= \sum\limits_{[i,j]}\sum\limits_{k} g_k^{[i,j]} = \sum\limits_{[i,j]}\sum\limits_{k} \left\|\boldsymbol{\mathcal{X}}^{[i,j]}_k\times_1 {\boldsymbol{A}_{[i,j]}}^* \times_2 {\boldsymbol{B}_k}^* \times_3 {\boldsymbol{C}^{[i,j]}_k}^*\right\|_F^2\\
    \text{s.t. }{\boldsymbol{A}_{[i,j]}}^*\boldsymbol{A}_{[i,j]} = \boldsymbol{I}_m, {\boldsymbol{B}_k}^* \boldsymbol{B}_k = \boldsymbol{I}_n, {\boldsymbol{C}^{[i,j]}_k}^*\boldsymbol{C}^{[i,j]}_k = \boldsymbol{I}_p,\text{ for all }i,j,k.
\end{array}.
\label{model: MT_max_model}
\end{equation}

For solving \eqref{model: MT_max_model}, similar to Algorithm \ref{alg: HOSVD Algorithm}, one can first use the singular vectors of matricization of $\boldsymbol{\mathcal{X}}^{[i,j]}_k$ in different mode as the initialization of $\boldsymbol{A}_{[i,j]}$, $\boldsymbol{B}_k$ and $\boldsymbol{C}^{[i,j]}_k$, respectively. Then update $\boldsymbol{A}_{[i,j]}$, $\boldsymbol{B}_k$ and $\boldsymbol{C}^{[i,j]}_k$  alternatively with the optimal solution of corresponding subproblem, e.g.,

\begin{itemize}
    \item For $\boldsymbol{A}_{[i,j]}$
    \begin{equation}\label{equ: update A}
\begin{aligned}
\boldsymbol{A}_{[i,j]} &= \argmax \limits_{{\boldsymbol{A}}^*\boldsymbol{A} = \boldsymbol{I}_m} \sum\limits_{k}\  g_k^{[i,j]}\left(\boldsymbol{A},\boldsymbol{B}_k,\boldsymbol{C}^{[i,j]}_k\right)\\
& = m \text{ principal eigenvectors of } \sum\limits_{k}\left(_1\boldsymbol{M}^{[i,j]}_k {_1\boldsymbol{M}^{[i,j]}_k}^*\right)\end{aligned}
\end{equation}
where $_1\boldsymbol{M}^{[i,j]}_k = \left(\boldsymbol{\mathcal{X}}^{[i,j]}_k\times_2 {\boldsymbol{B}_k}^* \times_3{\boldsymbol{C}^{[i,j]}_k}^*\right)_{(1)}$.
\item For $\boldsymbol{B}_k$
\begin{equation}\label{equ: update B}
\begin{aligned}
    \boldsymbol{B}_{k} &= \argmax\limits_{{\boldsymbol{B}}^*\boldsymbol{B} = \boldsymbol{I}_n}\ \sum\limits_{[i,j]}\ g_k^{[i,j]} \left(\boldsymbol{A}_{[i,j]},\boldsymbol{B},\boldsymbol{C}^{[i,j]}_k\right)\\
    & = n \text{ principle eigenvectors of } \sum\limits_{[i,j]} \left(_2\boldsymbol{M}_{k}^{[i,j]} {_2\boldsymbol{M}_{k}^{[i,j]}}^*\right)
\end{aligned}
\end{equation}
where $_2\boldsymbol{M}_{k}^{[i,j]} = \left(\boldsymbol{\mathcal{X}}^{[i,j]}_k\times_1 {\boldsymbol{A}_{[i,j]}}^* \times_3{\boldsymbol{C}^{[i,j]}_k}^* \right)_{(2)}$.
\item For $\boldsymbol{C}^{[i,j]}_k$
\begin{equation}\label{equ: update C}
\begin{aligned}
    \boldsymbol{C}^{[i,j]}_k &= \argmax\limits_{{\boldsymbol{C}}^*\boldsymbol{C} = \boldsymbol{I}_p}\ g_k^{[i,j]}\left(\boldsymbol{A}_{[i,j]},\boldsymbol{B}_k,\boldsymbol{C}\right)\\
    & = p \text{ principle eigenvectors of } _3\boldsymbol{M}^{[i,j]}_k {_3\boldsymbol{M}^{[i,j]}_k}^*
\end{aligned}
\end{equation}
where $_3\boldsymbol{M}_{k}^{[i,j]} = \left(\boldsymbol{\mathcal{X}}^{[i,j]}_k\times_1 {\boldsymbol{A}_{[i,j]}}^*\times_2 {\boldsymbol{B}_k}^*\right)_{(3)}$.
\end{itemize}

Details of updates are summarized in Algorithm\ref{alg: HOSVD_MT}.
\begin{algorithm}[H]
	\renewcommand{\algorithmicrequire}{{Input:}}
	\renewcommand{\algorithmicensure}{{Output:}}
	\caption{Groupwise Tucker compression algorithm}
	\label{alg: HOSVD_MT}
	\begin{algorithmic}[1]
	    \STATE Initialize $(\boldsymbol{A}_{[i,j]})_0$, $(\boldsymbol{B}_k)_0$ and $(\boldsymbol{C}^{[i,j]}_k)_0$ for $[i, j]\in [J]\times [K_i]$ and $k\in [J]$.
		\FOR{$s = 0,\dots, S-1$}
		\FOR{$[i,j] \in [J]\times [K_i]$}
        \STATE Compute $(\boldsymbol{A}_{[i,j]})_{s+1}$ through \eqref{equ: update A} with $(\boldsymbol{B}_k)_s$ and $(\boldsymbol{C}^{[i,j]}_k)_s$, for $k\in [J]$. 
        \ENDFOR
        \FOR{$k \in [J]$}
        \STATE Compute $(\boldsymbol{B}_k)_{s+1}$ through \eqref{equ: update B} with $(\boldsymbol{A}_{[i,j]})_{s+1}$ and $(\boldsymbol{C}^{[i,j]}_k)_s$, for $[i, j]\in [J]\times [K_i]$
        \ENDFOR
        \FOR{$(k,[i,j])\in [J]\times [J]\times [K_i]$}
        \STATE Compute $(\boldsymbol{C}^{[i,j]}_k)_{s+1}$ through \eqref{equ: update C} with  $(\boldsymbol{A}_{[i,j]})_{s+1}$ and $(\boldsymbol{B}_k)_{s+1}$.
        \ENDFOR
        \ENDFOR
        \FOR{$(k,[i,j])\in [J]\times [J]\times [K_i]$}
        \STATE $(\boldsymbol{\mathcal{G}}^{[i,j]}_k)_S = \boldsymbol{\mathcal{X}}^{[i,j]}_k\times_1 (\boldsymbol{A}_{[i,j]})_S^* \times_2 (\boldsymbol{B}_k)_S^* \times_3 (\boldsymbol{C}^{[i,j]}_k)_S^*$.
       \ENDFOR
        \ENSURE  $(\boldsymbol{A}_{[i,j]})_S,(\boldsymbol{B}_k)_S, (\boldsymbol{C}^{[i,j]}_k)_S$ and  $(\boldsymbol{\mathcal{G}}^{[i,j]}_k)_S$ for $[i, j]\in [J]\times [K_i]$ and $k\in [J]$.
    \end{algorithmic}  
\end{algorithm}

In Algorithm {\ref{alg: HOSVD_MT}}, it follows a Gauss-Seidel methodology that the factor matrices $\boldsymbol{A}_{[i,j]}$,$\boldsymbol{B}_k$ and $\boldsymbol{C}^{[i,j]}_k$ are alternatingly updated. The computation of the core tensors $\boldsymbol{\mathcal{G}}^{[i,j]}_k$ is not required in each iteration but is performed once after completion. The algorithm will converge to a solution where the objective function $f$ ceases to decrease. However, similar to the HOOI algorithm \cite{de2000best}, this method is not guaranteed to converge to the global optimum. Usually, one can use a predefined iteration number $S$ as a stopping criterion.

\subsection{Storage and SINR Complexity Reduction}{\label{section: SINR Complexity Analysis}}
For ease of exposition, we assume the numbers of users in all BSs are equal, denoted as $K$. Then with $J$ BSs there are $J^2 K$ channels in the entire system. So before compression, it is necessary to store $J^2 K$ complete channel tensors, each with a size of $M\times N\times P$. After compression using {(\ref{model: MT_origin})}, there are $JK$ first-mode factor matrices $\boldsymbol{A}_{[i,j]}$ with dimensions $M\times m$, $J$ second-mode factor matrices $\boldsymbol{B}_k$ with dimensions $N\times n$, $J^2 K$ third-mode factor matrices $\boldsymbol{C}^{[i,j]}_k$ with dimensions $P\times p$ and core tensors $\boldsymbol{\mathcal{G}}^{[i,j]}_k$ with size of $m\times n\times p$. Table \ref{tab: storage complexity} outlines the storage requirements before and after compression.
\begin{table}[H]
    \centering
    \begin{tabular}{|c|c|c|}
        \hline 
        & Original & Compressed  \\
        \hline
        Storage &  $J^2 K M N P$ & $J^2 K (m n p + P p) + J K M m + J N n$ \\
        \hline
    \end{tabular}
    \caption{ The storage (measured in complex double-precision floating-point numbers) of whole system before and after compression.}
    \label{tab: storage complexity}
\end{table}

For the SINR calculation, the first step is to reconstruct the channel matrices from channel tensors. Then is to compute the precoding matrix $\boldsymbol{P}^{[i,j]}$, covariance matrix $\boldsymbol{Q}^{[i,j]}$ and its inversion, along with the filter matrix $\boldsymbol{W}^{[i,j]}$. The SINR is finally determined using {(\ref{equ: SINR formula})}. The comparative computational complexity for each step before and after compression is outlined below.

\begin{enumerate}
    \item[$\bullet$] For reconstructing channel matrices $\boldsymbol{H}^{[i,j]}_k/ \boldsymbol{\widetilde{H}}^{[i,j]}_k$. According to {(\ref{equ: compressed matrix SINR reduction})}, before compression, directly performing a tensor-vector product $\boldsymbol{\mathcal{X}}^{[i,j]}_k\times_3 {\boldsymbol{c}_k^{[i,j]}}^*$ requires $MNP$ floating-point operations (flops) . After compression,  $(\boldsymbol{\mathcal{G}}^{[i,j]}_k\times_3 ({\boldsymbol{c}_k^{[i,j]}}^* \boldsymbol{C}^{[i,j]}_k))$ involves a matrix-vector product and a tensor-vector product, requiring $mnp + Pp$ flops.
    
    \item[$\bullet$] For the precoding matrix  $\boldsymbol{V}^{[k,l]}/ \boldsymbol{\widetilde{V}}^{[k,l]}$. Before compression, directly computing truncated-$L$ SVD of channel matrix $\boldsymbol{H}_{k}^{[k,l]}$ requires $O(MNL)$ flops.  After compression, only the truncated-$L$ SVD of  $\boldsymbol{\widetilde{H}}_{k}^{[k,l]}$ is needed due to the orthogonality of $\boldsymbol{A}_{[k,l]}$ and $\boldsymbol{B}_k$. It requires $O(mnL)$ flops to obtain $ \boldsymbol{\widetilde{V}}^{[k,l]}$, and the relationship between $\boldsymbol{V}^{[k,l]}$ and $\boldsymbol{\widetilde{V}}^{[k,l]}$ is: 
    \begin{equation}
     \boldsymbol{V}^{[k,l]} \approx \boldsymbol{B}_k \boldsymbol{\widetilde{V}}^{[k,l]},
    \end{equation}
    thereby eliminating the right factor matrix $\boldsymbol{B}_k$ in precoding process:
    \begin{equation}
     \boldsymbol{H}_{k}^{[i,j]} \boldsymbol{V}^{[k,l]}\approx \boldsymbol{A}_{[i,j]} \boldsymbol{\widetilde{H}}_{k}^{[i,j]} \boldsymbol{\widetilde{V}}^{[k,l]}
     \label{equ: compressed pre}
    \end{equation}
    \item[$\bullet$] For the covariance matrix $\boldsymbol{Q}^{[i,j]}/\boldsymbol{\widetilde{Q}}^{[i,j]}$ and its inversion.
    Plug   {(\ref{equ: compressed pre})} into   {(\ref{equ: covariance matrix})},
    \begin{equation}\label{equ: Cov}
    \begin{array}{rcl}
        \boldsymbol{Q}^{[i,j]} & = &  \sum\limits_{k = 1}^J \sum\limits_{l = 1}^{K} \boldsymbol{H}_{k}^{[i,j]} \boldsymbol{V}^{[k,l]}(\boldsymbol{H}_{k}^{[i,j]} \boldsymbol{V}^{[k,l]})^* + \sigma^2\boldsymbol{I}_M\\
         & \approx & \boldsymbol{A}_{[i,j]}(\sum\limits_{k = 1}^J \sum\limits_{l = 1}^{K} \boldsymbol{\widetilde{H}}_{k}^{[i,j]} \boldsymbol{\widetilde{V}}^{[k,l]}(\boldsymbol{\widetilde{H}}_{k}^{[i,j]} \boldsymbol{\widetilde{V}}^{[k,l]})^*) {\boldsymbol{A}_{[i,j]}}^* + \sigma^2 \boldsymbol{I}_M,
    \end{array}
    \end{equation}
    Apply Woodbury formula \cite{yip1986note} to the last equation of {(\ref{equ: Cov})}, the inversion $(\boldsymbol{Q}^{[i,j]})^{-1}$ is given by 
    \begin{equation}
        (\boldsymbol{Q}^{[i,j]})^{-1} \approx \frac{1}{\sigma^2} \boldsymbol{I}_M - \frac{1}{\sigma^2} \boldsymbol{A}_{[i,j]} (\boldsymbol{\widetilde{Q}}^{[i,j]})^{-1}(  \boldsymbol{\widetilde{Q}}^{[i,j]} - \sigma^2 \boldsymbol{I}_m) {\boldsymbol{A}_{[i,j]}}^*,
        \label{equ: reduced inversion of interference matrix}
    \end{equation}
    where
    \begin{equation*}
        \boldsymbol{\widetilde{Q}}^{[i,j]} = \sum\limits_{k = 1}^J \sum\limits_{l = 1}^{K} \boldsymbol{\widetilde{H}}_{k}^{[i,j]} \boldsymbol{\widetilde{V}}^{[k,l]} (\boldsymbol{\widetilde{H}}_{k}^{[i,j]} \boldsymbol{\widetilde{V}}^{[k,l]})^* + \sigma^2 \boldsymbol{I}_m.
    \end{equation*}
    here the inversion of the $\boldsymbol{Q}^{[i,j]}$ is substituted with the inversion of a smaller matrix $ \boldsymbol{\widetilde{Q}}^{[i,j]}$ and several matrix-matrix multiplications. The matrix $\boldsymbol{A}_{[i,j]}$ does not need to be explicitly multiplied as it can be eliminated in subsequent computations owing to its orthogonality.  
    
    In {(\ref{equ: Cov})}, the  matrix-matrix product $\boldsymbol{H}_{k}^{[i,j]} \boldsymbol{V}^{[k,l]}$ and $\boldsymbol{\widetilde{H}}_{k}^{[i,j]} \boldsymbol{\widetilde{V}}^{[k,l]}$ requires $MNL$ flops and $mnL$ flops, respectively. Then $\boldsymbol{H}_{k}^{[i,j]} \boldsymbol{V}^{[k,l]}(\boldsymbol{H}_{k}^{[i,j]} \boldsymbol{V}^{[k,l]})^*$ and $\boldsymbol{\widetilde{H}}_{k}^{[i,j]} \boldsymbol{\widetilde{V}}^{[k,l]} (\boldsymbol{\widetilde{H}}_{k}^{[i,j]} \boldsymbol{\widetilde{V}}^{[k,l]})^*$ requires $M^2 L$ flops and $m^2 L$ flops, respectively. Consequently,  calculating $\boldsymbol{Q}^{[i,j]}$ and $\boldsymbol{\widetilde{Q}}^{[i,j]}$ requires $J K(M N L + M^2 L)$ and $J K (m n L + m^2 L)$ flops, respectively. $(\boldsymbol{Q}^{[i,j]})^{-1}$ and $\boldsymbol{\widetilde{Q}}^{[i,j]}$ requires $O(M^3)$ and $O(m^3)$ flops, respectively. The matrix-matrix multiplications $(\boldsymbol{\widetilde{Q}}^{[i,j]})^{-1}(  \boldsymbol{\widetilde{Q}}^{[i,j]} - \sigma^2 \boldsymbol{I}_m)$ requires $m^3$ flops.
    \item[$\bullet$] For the filter matrix $\boldsymbol{W}^{[i,j]}/\boldsymbol{\widetilde{W}}^{[i,j]}$.
    By substituting   {(\ref{equ: reduced inversion of interference matrix})} into {(\ref{equ: MMSE receive filter})}, $\boldsymbol{W}^{[i,j]}$ can be expressed as:
    \begin{equation}
        \boldsymbol{W}^{[i,j]} \approx  \boldsymbol{A}_{[i,j]} \boldsymbol{\widetilde{W}}^{[i,j]},
    \end{equation}
    where $\boldsymbol{\widetilde{W}}^{[i,j]} = \frac{1}{\sigma^2}(\boldsymbol{I}_m - (\boldsymbol{\widetilde{Q}}^{[i,j]})^{-1} (\boldsymbol{\widetilde{Q}}^{[i,j]} - \sigma^2 \boldsymbol{I}_m)\boldsymbol{\widetilde{H}}_{i}^{[i,j]} \boldsymbol{\widetilde{V}}^{[i,j]})$.
    In this procedure, the matrix is still retained  $\boldsymbol{A}_{[i,j]}$ and not explicitly multiplied. The matrix $(\boldsymbol{\widetilde{Q}}^{[i,j]})^{-1} (\boldsymbol{\widetilde{Q}}^{[i,j]} - \sigma^2 \boldsymbol{I}_m)$ and $\boldsymbol{\widetilde{H}}_{i}^{[i,j]} \boldsymbol{\widetilde{V}}^{[i,j]}$ have already been obtained in the previous procedure. Consequently, the computational complexity for $\boldsymbol{W}^{[i,j]}$ and $\boldsymbol{\widetilde{W}}^{[i,j]}$ is $M^2 L$ and $m^2 L + mL$ flops, respectively.
    \item[$\bullet$] For $\boldsymbol{s}^{[i,j]}_r/ \boldsymbol{\widetilde{s}}^{[i,j]}_r$ and ${SINR}_r^{[i,j]}/ \widetilde{SINR}_r^{[i,j]}$. We can finally eliminate the factor matrix $\boldsymbol{A}_{[i,j]}$. 
    The $\boldsymbol{s}^{[i,j]}_r$ in (\ref{equ: SINR formula}) only involves a vector inner product, thus it requires $M$ flops while for $\boldsymbol{\widetilde{s}}^{[i,j]}_r$,
    \begin{equation}
    \begin{aligned}
        \boldsymbol{\widetilde{s}}^{[i,j]}_r & = \left| \left[\boldsymbol{\widetilde{W}}^{[i,j]}\right]_r^* {\boldsymbol{A}_{[i,j]}}^* \boldsymbol{A}_{[i,j]} \left[\boldsymbol{\widetilde{H}}^{[i,j]}_i \boldsymbol{\widetilde{V}}^{[i,j]}\right]_r \right|^2\\
        & = \left| \left[\boldsymbol{\widetilde{W}}^{[i,j]}\right]_r^* \left[\boldsymbol{\widetilde{H}}^{[i,j]}_i \boldsymbol{\widetilde{V}}^{[i,j]}\right]_r \right|^2
    \end{aligned}
    \end{equation}
    it requires $m$ flops. 
    For ${SINR}_r^{[i,j]}$, the $\left[\boldsymbol{W}^{[i,j]}\right]_r^* \boldsymbol{Q}^{[i,j]}\left[\boldsymbol{W}^{[i,j]}\right]_r$ in \eqref{equ: SINR formula} needs a matrix-vector product and a vector-vector product, which needs $M^2 + M$ flops. Similarly, $\widetilde{SINR}_r^{[i,j]}$ needs $m^2 + m$ flops.
\end{enumerate}

\begin{figure}[H]
    \centering
    \includegraphics[width = 0.6\textwidth]{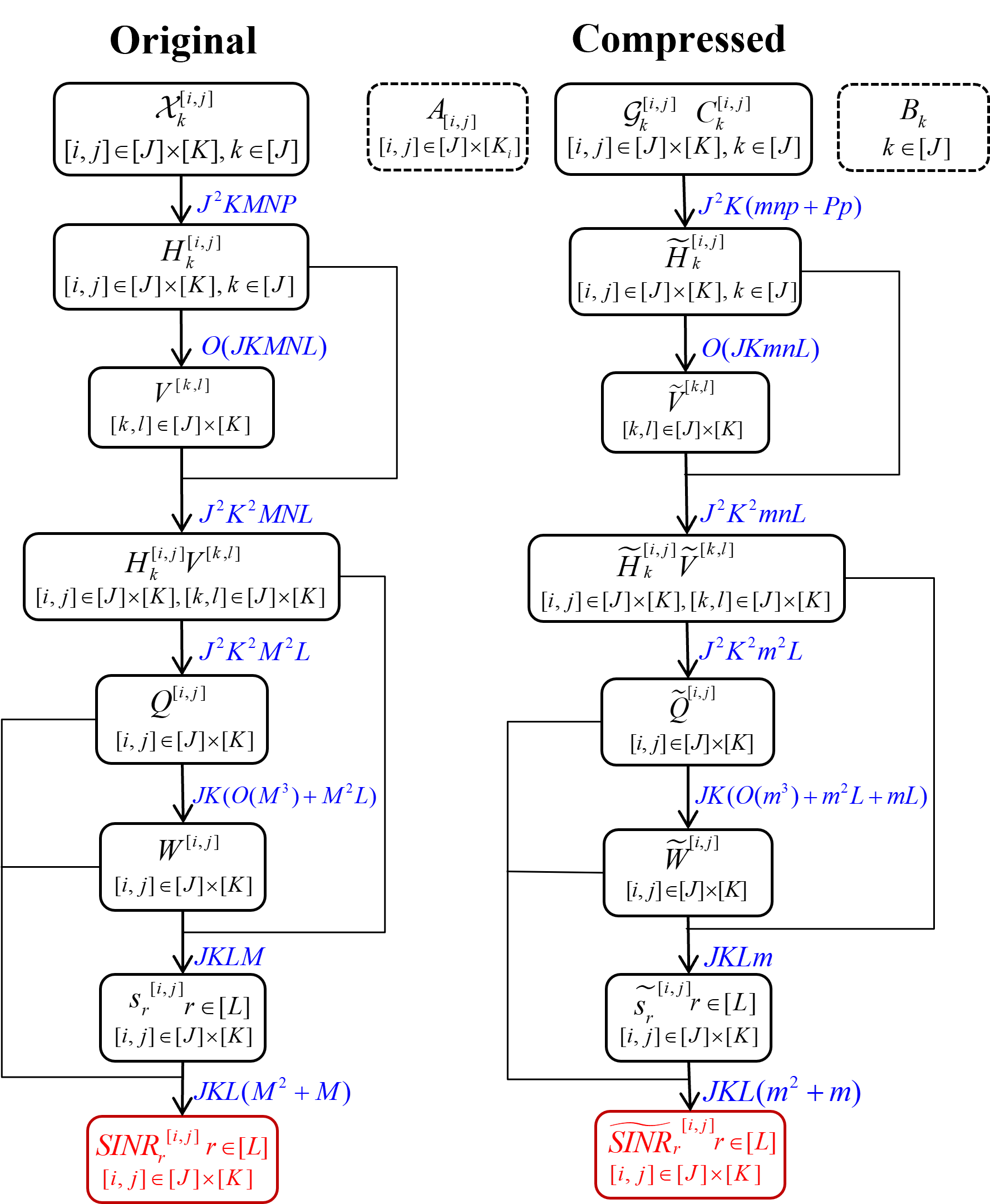}
    \caption{Computational complexity comparison.}
    \label{fig: Computational complexity}
\end{figure}

Clearly, the SINR calculation in {(\ref{equ: SINR formula})} can be directly evaluated using the compressed channel matrix $\boldsymbol{\widetilde{H}}_{k}^{[i,j]}$.
A comparison of the total computational complexity for SINR using complete and  compressed channel matrices is presented in Figure \ref{fig: Computational complexity}.

Indeed, utilizing the compressed channel matrix $\boldsymbol{\widetilde{H}}^{[i,j]}_k$ with a smaller size leads to a reduction in computational complexity. However, this decrease in size comes at the cost of increased approximation error. Consequently, choosing the compressed channel size necessitates careful consideration to strike an optimal balance between computational efficiency and accuracy. The subsequent section will delve into this trade-off through a series of experiments.

\section{Numerical Experiments}{\label{section: Numerical Results}}

\begin{figure}[h!]
    \centering
    \includegraphics[width = .45\linewidth]{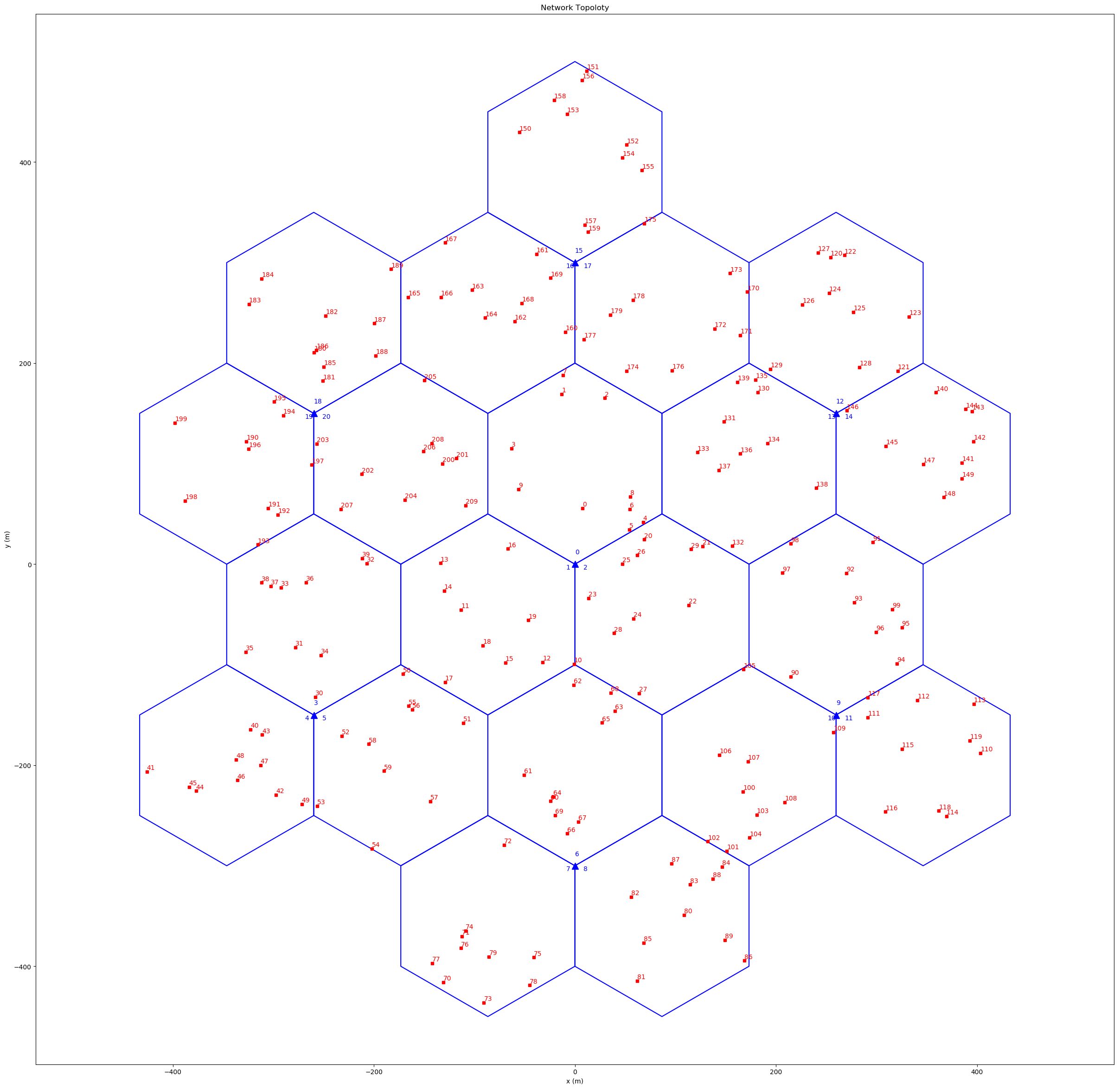}
    \caption{ Topological diagram for 21 BSs and 210 UEs MIMO system.}
    \label{fig: topo MIMO system}
\end{figure}

In this section, we conduct experiments on MIMO systems of different scales and channels to evaluate the numerical performance of the proposed compression models. The channel matrices used in the experiments are generated based on the specifications provided by the 3rd Generation Partnership Project (3GPP) \cite{3gpp.38.901} at a fixed time point. We consider two popular channel matrix sizes: $64\times512$ (64 receiving and 512 transmitting antennas) and $8\times256$, the submatrices number $P$ is $401$. Figure {\ref{fig: topo MIMO system}} provides a pictorial depiction of a MIMO system comprising 21 base stations and 10 users per cell used in our experiments. In the experiments, we simplify the scenario by considering only the inter-cell interference from the first user equipment (UE) of other cells and ignoring the intra-cell interference.

To evaluate the performance of the compression models in terms of reducing storage requirements and accelerating SINR computation, we introduce two compression ratios: $R_s$ and $R_t$. The compression ratio $R_s$ is defined as the ratio of memory usage for channel data before and after compression (as shown in Table {\ref{tab: storage complexity}}). $R_s$ can be directly calculated by the following equation.
\begin{equation}\label{equ: compression ratio}
R_s = \frac{J^2 K M N P}{J^2 K(m n p + P p) + J K M m + J N n}
\end{equation}

For the speedup ratio $R_t$, we count the run time of calculating SINR as depicted in Figure \ref{fig: Computational complexity}. Assume $t_1$ and $t_2$ are the run time for original channel data and compressed data, respectively, then $R_t = t_1/t_2$.

We evaluate the accuracy of the compression models by calculating the mean relative error $e_c$ of the SINR values before and after compression. The mean relative error is computed as the average of the absolute differences between the original SINR values $s_r^{[i,j]}$ and the compressed SINR values $\widetilde{s}_r^{[i,j]}$ for all $r\in [L], [i,j]\in [J]\times[K]$, divided by the original SINR values, e.g.,

\begin{equation}\label{equ: SINR error}
e_c = \text{mean of }\left\{\left|\frac{{SINR}_r^{[i,j]} - \widetilde{SINR}_r^{[i,j]}}{{SINR}_r^{[i,j]}}\right|\right\}.
\end{equation}

From Tables {\ref{tab: storage complexity}} and Figure \ref{fig: Computational complexity}, we can observe a strong correlation between $R_t$ and $R_s$. Smaller compressed channel sizes lead to higher $R_s$ and $R_t$ but lower accuracy due to increased approximation errors. In our experiments, we compare three compression models: Individual Tucker model in {(\ref{model: Tucker decomposition})}, Shared Tucker model in {(\ref{model: MT1})} and Groupwise Tucker  model in {(\ref{model: MT_origin})}. We firstly maintain the mean relative error $e_c$ at around $10\%$ and adjust the compression parameters $m, n, p$ to achieve higher values of $R_s$ and $R_t$ in the Groupwise Tucker model. Since $M$ is already small, we focus on tuning the parameters $n$ and $p$ while setting $m$ to be equal (or close) to $M$. Subsequently, we set the same compression parameters $m, n, p$ in the Individual Tucker model and Shared Tucker model. The results of three compression models are summarized in Table {\ref{tab: Compression Results for 21 base stations 210 users data.}}.
\begin{table}[H]
\centering
\begin{tabular}{ccccclccc}
\toprule
\multicolumn{5}{c}{Settings}                                                                                                         &  & \multicolumn{3}{c}{Results}        \\ \cline{1-5} \cline{7-9}
$L$                   & $(M,N,P)$                         & $K$                   & $(m,n,p)$                         & Model        &  & $R_t$     & $R_s$    & $e_c$       \\ \cline{1-5} \cline{7-9} 
\multirow{6}{*}{$21$} & \multirow{6}{*}{$(64, 512, 401)$} & \multirow{3}{*}{$5$}  & \multirow{3}{*}{$(60, 230, 150)$} & Individual Tucker   &  & $0.17951$ & $5.8354$ & $5.6042\%$  \\
                      &                                   &                       &                                   & Shared Tucker &  & $5.7931$  & $6.1684$ & $12.5947\%$ \\
                      &                                   &                       &                                   & \textbf{Groupwise Tucker}   &  & $6.1904$  & $6.1648$ & $9.3929\%$  \\ \cline{3-5} \cline{7-9} 
                      &                                   & \multirow{3}{*}{$10$} & \multirow{3}{*}{$(60, 270, 190)$} & Individual Tucker       &  & $0.17072$ & $3.9863$ & $2.9760\%$  \\
                      &                                   &                       &                                   & Shared Tucker     &  & $3.3805$  & $4.1658$ & $13.2991\%$ \\
                      &                                   &                       &                                   & \textbf{Groupwise Tucker}   &  & $4.1080$  & $4.1648$ & $8.2277\%$  \\ \cline{1-5} \cline{7-9} 
\multirow{6}{*}{$21$} & \multirow{6}{*}{$(8, 256, 401)$}  & \multirow{3}{*}{$5$}  & \multirow{3}{*}{$(8, 130, 120)$}  & Individual Tucker       &  & $0.15937$ & $3.9815$ & $6.5190\%$  \\
                      &                                   &                       &                                   & Shared Tucker     &  & $1.2023$  & $4.7489$ & $13.3681\%$ \\
                      &                                   &                       &                                   & \textbf{Groupwise Tucker }  &  & $1.3500$  & $4.7405$ & $10.8535\%$ \\ \cline{3-5} \cline{7-9} 
                      &                                   & \multirow{3}{*}{$10$} & \multirow{3}{*}{$(8, 140, 140)$}  & Individual Tucker       &  & $0.14927$ & $3.3003$ & $3.9147\%$  \\
                      &                                   &                       &                                   & Shared Tucker     &  & $1.1768$  & $3.8566$ & $9.9698\%$  \\
                      &                                   &                       &                                   & \textbf{Groupwise Tucker}   &  & $1.1014$  & $3.8536$ & $9.0037\%$  \\ \toprule 
\end{tabular}
\caption{ Compression results of Groupwise Tucker , Individual Tucker and Shared Tucker models for MIMO systems with different scales and channels.}
\label{tab: Compression Results for 21 base stations 210 users data.}
\end{table}

In rows 3, 6 (or rows 9, 12) of Table {\ref{tab: Compression Results for 21 base stations 210 users data.}}, we fix $M, N, P$ and consider three different numbers of UEs $K$ in each BS. We tune $m, n, p$ such that the relative error $e_c$ is below $10\%$. The results indicate that larger compression sizes $m, n, p$ (resulting in smaller $R_t$ and $R_s$) are required to ensure the relative error $e_c$ when more users are considered.

In rows 1-3 (or rows 4-6, rows 7-9, rows 10-12) of Table {\ref{tab: Compression Results for 21 base stations 210 users data.}}, we compare the three compression models under the same settings. 
The error of the Individual Tucker decomposition is smaller than that of the Groupwise Tucker model and Shared Tucker model, but the speedup ratio $R_t$ of the Individual Tucker decomposition is less than 1, indicating that it requires more computational cost compared to the original data. Comparing the Groupwise Tucker and Shared Tucker models, they have similar performance on accelerating the SINR process with $R_t$ greater than 1, while the solution of Shared Tucker model is more accurate as the relative error $e_c$ is notably smaller than that of the Groupwise Tucker model. 

For the compression ratio $R_s$, the Shared Tucker model has a slightly higher value than the Groupwise Tucker model. This is consistent with their factorization forms in their respective models. Comparing rows 1-6 with rows 7-12 of Table {\ref{tab: Compression Results for 21 base stations 210 users data.}} respectively, we can observe that the compression performance is better when the channel size is larger, particularly in terms of the speedup ratio $R_t$.

\begin{figure}[H]
    \centering
    \includegraphics[width = .65\linewidth]{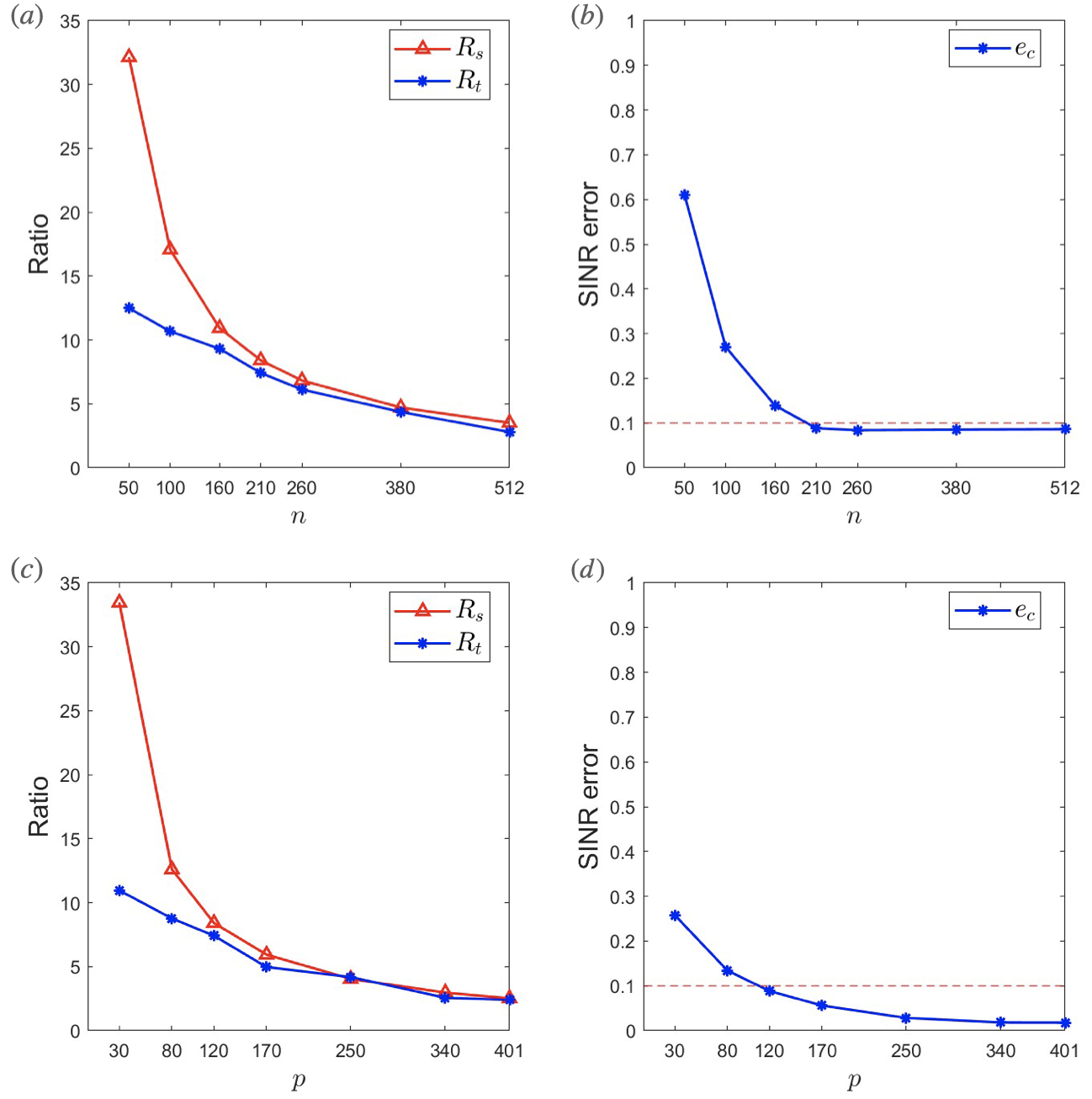}
    \caption{The compression results for channel data with $L = 21$, $K = 2 $ and $(M, N, P) = (64, 512, 401)$. The red dashed lines in Figure (b) and Figure (d) represent the $10\%$ mean relative error we desire.}
    \label{fig: change_Col_Num}
\end{figure}
In addition to the previous results, we also provide the univariate tuning results of Groupwise Tucker  compression model to show the impact of changing the number of compressed columns (or submatrices) on the performance metrics. We use the channel data with $L=21$ and $K=2$ and the antenna configuration is $(M,N,P)=(64,512,401)$.
We first fix $(m, p) = (60, 120)$ and vary the compressed columns number $n$, the result for $(R_t, R_c)$ and $e_c$ are shown in Figure \ref{fig: change_Col_Num}(a) and Figure \ref{fig: change_Col_Num}(b), respectively. Also we consider to fix $(m, n)$
and vary $p$, the corresponding results are depicted in Figure {\ref{fig: change_Col_Num}}(c) and Figure {\ref{fig: change_Col_Num}}(d).

Upon inspection of {(\ref{equ: compression ratio})} derived from the storage complexity in Table {\ref{tab: storage complexity}}, we can deduce that $R_s$ is inversely proportional to $p$ and $n$. The numerical observations depicted in (a) and (c) of Figure {\ref{fig: change_Col_Num}} align with our theoretical predictions. Regarding the SINR error plots in Figure {\ref{fig: change_Col_Num}} (b) and (d), it appears that smaller values of $n$ or $p$ are sufficient to achieve the best SINR error performance. This suggests that a lower number of compressed columns or submatrices can effectively preserve the accuracy of the SINR calculation.

\section{Conclusion and Discussion}{\label{section: Conclusion and Discussion}}
In this paper, we propose the Groupwise Tucker  compression model based on Tucker decomposition to compress channel data in massive MIMO systems. Our approach enables faster signal transmission operations, particularly in SINR calculation, while significantly reducing storage requirements. The numerical results demonstrate substantial improvements in both storage efficiency and SINR calculation speed with high precision. However, there are two remaining issues that require further investigation. 

Firstly, the current algorithm relies on SVD, which can be time-consuming due to the high dimensions of channel data in our problem. Although the compression task can be performed offline, reducing this computational cost is worth pondering. One possible solution is to incorporate SVD-free techniques into the framework, as suggested in a recent paper \cite{xiao2021efficient}.

Secondly, the selection of compression parameters currently involves a time-consuming trial-and-error approach. While we have explored using evolutionary algorithms such as Particle Swarm Optimization (PSO) \cite{kennedy1995particle} and surrogate learning to address this issue, none of these methods have yielded satisfactory results. Therefore, another avenue for improvement is to develop rank-selection approaches for determining compression parameters more efficiently.

Addressing these two issues would further enhance the practicality and efficiency of the Groupwise Tucker  compression model, making it more viable for real-world applications.

\bibliographystyle{unsrt}
\bibliography{Template}

\end{document}